\documentclass[a4paper,12pt]{article}

\usepackage{latexsym}
\usepackage{amsfonts}
\newtheorem{thm}{Theorem}[section]
\newtheorem{exa}[thm]{Example}

\newtheorem{cor}[thm]{Corollary}

\newcommand{\qed}{\unskip\protect\nolinebreak\mbox{\quad$\Box$}\vspace{3mm}}
 
\newcommand{\Z}{\mathbb{Z}}

\newcommand{\psl}{\mathop{\rm PSL}\nolimits}
\newcommand{\pgl}{\mathop{\rm PGL}\nolimits}

\makeatletter
\def\imod#1{\allowbreak\mkern10mu({\operator@font mod}\,\,#1)}
\makeatother

\usepackage{fullpage}

\begin{document}

\title{Terraces for small groups} 
\author{M. A. Ollis\footnote{E-mail address: \texttt{matt@marlboro.edu}.}  \\ \\
Marlboro College, P.O.~Box A, Marlboro,  \\ Vermont 05344, USA}

\date{}

\maketitle

\begin{abstract}
We use heuristic algorithms to find terraces for small groups.  We show that Bailey's Conjecture (that all groups other than the non-cyclic elementary abelian 2-groups are terraced) holds up to order~511, except possibly at orders 256 and 384.    We also show that Keedwell's Conjecture (that all non-abelian groups of order at least~10 are sequenceable) holds up to order~255, and for the groups $A_6$, $S_6$,  $\psl(2,q_1)$ and $\pgl(2, q_2)$ where $q_1$ and~$q_2$ are prime powers with $3 \leq q_1 \leq  11$ and $3 \leq q_2 \leq 8$.  A sequencing for a group of a given order implies the existence of a complete latin square at that order.  We show that there is a sequenceable group for each odd order up to~555 at which there is a non-abelian group.  This gives 31 new orders at which complete latin squares are now known to exist, the smallest of which is 63.  In addition, we consider terraces with some special properties, including constructing a directed $T_2$-terrace for the non-abelian group of order~21 and hence a Roman-2 square of order~21 (the first known such square of odd order).  Finally we report  the total number terraces and directed terraces for groups of order at most~15.

\vspace{3mm}


\end{abstract}

\section{Introduction}

Let $G$ be a group of order~$n$.  Let $(a_1, a_2, \ldots, a_n)$ be an arrangement of the elements of $G$ and define ${\bf b} = (b_1, b_2, \ldots, b_{n-1})$ by $b_i = a_i^{-1}a_{i+1}$ for each $i$ with $1 \leq i \leq n-1$.  If the elements of ${\bf b}$ are distinct then ${\bf a}$ is a {\em directed terrace} for $G$ and ${\bf b}$ is a {\em sequencing}.  A group that has a directed terrace and sequencing is called {\em sequenceable}.  

Gordon \cite{Gordon61} introduced the concept as a tool for constructing ``complete Latin squares";  the same ideas had been used earlier in cyclic groups \cite{Williams49}.  A Latin square of order $n$ is an $n \times n$ array of $n$ symbols with each symbol appearing once in each row and once in each column.  A Latin square is {\em row-complete}, and is also known as a {\em Roman square}, if each ordered pair of symbols appears in consecutive positions within rows exactly once.  If the transpose is also row-complete, the square is {\em complete}. Gordon showed that if $( a_1, a_2, \ldots, a_n)$ is a directed terrace for a group of order $n$ then the Latin square with $(i,j)$ entry $a_i^{-1}a_j$ is complete.

A {\em binary group} is defined to be a group with a single involution. An abelian group is sequenceable if and only if it is a binary group and the three non-abelian groups of orders~6 and~8 are not sequenceable \cite{Gordon61}.  Many non-abelian groups and families of non-abelian groups have since been shown to have directed terraces, for example the dihedral and dicyclic groups of order at least~10.
No further exceptions are known.  See \cite{survey} for a survey of the current state of affairs.  

{\em Keedwell's Conjecture} is that all non-abelian groups of order at least~10 are sequenceable.  In the next section we show that it holds up to order~255 and that for each odd order of at most~555 at which there is a non-abelian group there is at least one sequenceable group.

Bailey \cite{Bailey84} generalised the notion of a sequencing (and introduced the terrace terminology) for use in constructing ``quasi-complete Latin squares".   A Latin square is {\em row quasi-complete} if each pair of symbols appears twice (in either order) in adjacent cells within rows of the square.  If the transpose of a row quasi-complete square is also row quasi-complete then the square is {\em quasi-complete}.   

With the arrangement ${\bf a}$ and associated list ${\bf b}$ as above, if each involution occurs exactly once in ${\bf b}$ and if for each $x \in G$ with $x^2 \neq e$ the sequence ${\bf b}$ contains $x$ and $x^{-1}$ twice in total then ${\bf a}$ is a {\em terrace} for $G$ and ${\bf b}$ is its associated {\em 2-sequencing}.   If $( a_1, a_2, \ldots, a_n)$ is a  terrace for a group of order $n$ then the Latin square with $(i,j)$ entry $a_i^{-1}a_j$ is quasi-complete \cite{Bailey84}.  

Non-cyclic elementary abelian 2-groups are not terraced \cite{Bailey84} and {\em Bailey's Conjecture} is that these are the only groups that do not have a terrace.  As with the directed case, many families of groups are known to be terraced, see \cite{survey}.  For example, Bailey's Conjecture holds for abelian groups \cite{OW4}.  In the next section we show that the conjecture is true up to order~511, with the possible exceptions of~256 and~384, and this and the results of Section~\ref{sec:special} imply the existence of terraces for new infinite families of groups. 

A sequence remains a terrace (or directed terrace) if we left-multiply all of its elements by any fixed group element.  Taking the element to be $a_1^{-1}$ gives a terrace that starts with $e$; such terraces are called {\em basic}.  Two terraces are {\em essentially different} if, after making them both basic by left-multiplying by a group element,  there is not an automorphism of the group that sends one to the other.

\begin{exa}\label{ex:lww}
Let~$\Z_n$ denote the additively written cyclic group of order~$n$ on the symbol set $\{0, 1, \ldots, n-1\}$.  The sequence
$$( 0, 1, n-1, 2, n-2, 3, n-3 \ldots )$$
is a terrace for~$\Z_n$.  When~$n$ is even it is directed.  The earliest publications of this construction seem to be~\cite{Lucas92} for even~$n$, where credit is given to Walecki, and~\cite{Williams49} for odd~$n$.
\end{exa}

In \cite{loseqs}, Anderson introduces a hill-climbing algorithm for finding directed terraces and he uses a variant of it in \cite{Anderson92} to look for terraces.  We describe our similar algorithms in the next section and report on which groups for which we were able to find terraces and sequencings.  The algorithms are implemented in GAP~\cite{gap} and are available at the following website:
\begin{center}
\texttt{http://cs.marlboro.edu/courses/matt/terraces}
\end{center}

In Section~\ref{sec:special} we consider terraces that have special properties that let them be used for general constructions of terraces for larger groups, for constructing combinatorial objects  or are of interest for their own sake.  Finally, in Section~\ref{sec:vsg} we report on exhaustive searches that take our full knowledge of how many essentially different terraces and directed terraces there are for all groups up to order~15.

Some groups we use: $D_{2m}$ is the dihedral group of order~$2m$; $Q_{4m}$ is the dicyclic group of order~$4m$, which is the quaternion group when~$m=2$ and the generalised quaternion group when~$m$ is a larger power of~2; $A_n$ and~$S_n$ are alternating and symmetric groups respectively; and $\psl(n,q)$ and~$\pgl(n,q)$, where~$q$ is a prime power, are the projective special and projective linear groups respectively.  

For other groups when there is no commonly used name we use the notation~$G_{n,i}$ for the group of order~$n$ in position~$i$ of GAP's Small Group Library \cite{gap}.  We use two such groups in examples and so give presentations for them here.  The first is the unique-up-to-isomorphism non-abelian group of order~21:
$$G_{21,1} = \langle u,v : u^7 = e = v^3,  vu = u^4v \rangle.$$
This is the smallest non-abelian group of odd order.   The second is one of the two non-abelian groups of order~27:
$$G_{27,4} = \langle u,v : u^9 = e = v^3,  vu = u^7v \rangle.$$

The proofs consist almost entirely of reports that the appropriate computer programs have run successfully.  As such they are omitted except for occasional ``proof notes" that contain some commentary.

\section{Keedwell's and Bailey's Conjectures}

We first consider directed terraces and Keedwell's Conjecture.   

For a hill-climbing algorithm we need three ingredients:  a search space through which we shall move, a neighbourhood that can be constructed for any element of the search space to give the potential next steps of our movement and a measure of our altitude (including knowing what altitude is sufficiently ``high" to be a solution to the problem).  We also include a ``teleport" function that moves us a significant distance across the search space without too much loss of altitude as a means for escaping local maxima.

For an arrangement ${\bf a} = (a_1, a_2, \ldots, a_n)$ of the elements of a group~$G$ of order~$n$, define ${\bf b}= (b_1, b_2, \ldots, b_{n-1})$ as usual by $b_i = a_i^{-1}a_{i+1}$ for each $i$ with $1 \leq i \leq n-1$.   The search space is the set of all arrangements of the group elements and the altitude of such an arrangement ${\bf a}$ is the number of distinct elements of ${\bf b}$.  The sequence ${\bf a}$ is a directed terrace if and only if its altitude is~$n-1$.

For movement we cut the sequence into a small number of pieces  and reassemble them in any order.  For example, if we cut $(a_1, a_2, \ldots a_n)$ at position $i$ and change the order of the two subsequences we reach the sequence 
$$(a_{i+1}, a_{i+2}, \ldots, a_n, a_1, a_2, \ldots a_i).$$
With two cuts there are five possible places to move that are different to our starting sequence.

These moves correspond well with our altitude measure as the portions of ${\bf b}$ generated by the subsequences do not change.  Using a move with $k$ cuts, the altitude can  be reduced at most $k$ and increased by up to $2k$.  This means we have a fairly smooth landscape.

The teleport function takes a random element of the sequence and moves it to the end.  This reduces the altitude by at most~2.  

The algorithm is now straightforward to describe.  
\begin{enumerate}
\item  Begin with a random arrangement of the elements of $G$.  
\item If the altitude is  $n-1$ then we have a directed terrace; return it.  Otherwise continue to Step 3.
\item Look through the neighbours of the arrangement one at a time and as soon as we see a higher one take that to be our new position and return to Step 2.  If there are no higher neighbours, continue to Step 4.
\item Teleport and return to Step 2.
\end{enumerate}
The number of cuts was limited to 2.  Experiments with 3 cuts suggested that the reduction in the number of steps to a solution was outweighed by the time taken to search for a beneficial 3-cut move.

Here is the collected outcome of the searches.  Anderson \cite{loseqs, nonab32} has already shown that it is true up to order~32 and for~$A_5$ and~$S_5$ using a similar algorithm.

\begin{thm}\label{th:kc}
All non-abelian groups of order $n$ with $10 \leq n \leq 255$ are sequenceable.  Also, the groups $A_6$, $S_6$,  $\psl(2,q_1)$ and $\pgl(2, q_2)$ are sequenceable where $q_1$ and~$q_2$ are prime powers with $3 \leq q_1 \leq  11$ and $3 \leq q_2 \leq 8$.
\end{thm}

As a directed terrace for a group of order~$n$ can be used to construct a complete Latin square of order~$n$ and there are many odd orders for which these are unknown, the program was pushed further on groups of odd order.  

\begin{thm}\label{th:odd}
For all odd orders~$n \leq 555$ for which there is a non-abelian group of order~$n$ there is at least one sequenceable group of order~$n$. 
\end{thm}

\noindent
Proof note.   Although theoretical constructions are known for several odd orders \cite{Keedwell81,Ollis14, Wang} the program was run and found directed terraces at these orders in any case.

At some values of~$n$ there are several non-abelian groups.  In these instances we found a sequencing for the first one in the ordering given by GAP's Small Group Library.
\qed

We collect the orders new to this work in the following result:

\begin{cor}\label{cor:newcls}
There is at least one sequenceable group at each of the orders
$$\begin{array}{l}
63, 105, 117, 135, 165, 171, 189, 195, 225, 231, 273, 275, 279, 285, 297, 315, \\
333, 351, 357, 385, 387, 405, 429, 459, 465, 483, 495, 513, 525, 549, 555.
\end{array}
$$
Hence there is a complete latin square of each of these orders.
\end{cor}

For Bailey's Conjecture, where the sequences have more freedom, we use more-or-less the same algorithm and are able to push to higher orders.

The changes to the algorithm are to the definition of altitude and the possible movement.

Given ${\bf a}$ and ${\bf b}$ with the usual notation, the altitude of ${\bf a}$ is now defined to be the sum of number of distinct involutions in ${\bf b}$ and the number of occurrences of  each element $g \in G$ with $g^2 \neq e$ up to a maximum of two from $g$ and $g^{-1}$ combined.   The altitude is $n-1$ if and only if ${\bf a}$ is a terrace. 

The idea of cutting and regluing remains the key to movement, but we now take advantage of the fact that if we reverse a subsequence of ${\bf a}$ then the elements of the portion of ${\bf b}$ that is generated are switched for their inverses.  This does not affect the altitude (for not-necessarily-directed terraces) and so after making one or two cuts as before we allow reversing the subsequences before re-assembly of the pieces.  For example, we might make one cut, at position~$i$, and reverse the second portion and glue the subsequences in the original order, to get:
$$(a_{1}, a_{2}, \ldots, a_i, a_n, a_{n-1}, \ldots a_{i+1}).$$
As in the directed case, making $k$ cuts can make the altitude at most~$k$ lower or at most~$2k$ higher, and we use the same teleport function that reduces the altitude by at most~2.  The algorithm is now essentially the same as the one given earlier.

Although constructions of directed terraces are available for infinitely many groups, they are not sufficiently dense to be worth factoring into the programs.  However, for terraces, we have the following very powerful result:

\begin{thm}\label{th:crcr} {\rm \cite{AandI92, AandI93}}
Let $G$ be a group and let $N$ be a normal subgroup of $G$. If $N$ has odd order and $G/N$ is terraced or if $N$ has odd index and $N$ is terraced, then $G$ is terraced.  
\end{thm}

This implies, for instance, that all groups of odd order are terraced.   Unfortunately, it does not help with 2-groups, which are by far the most abundant groups in the range we are considering.  However, among the 33,080 non-abelian non-2-groups of order at most~511 there are just 782 that cannot be terraced by using Theorem~\ref{th:crcr} provided Bailey's Conjecture holds at all proper divisors of the order.  Further, 29 of these 782 are dihedral groups of order~$4p$, for $p$ an odd prime, and these groups are known to be terraced \cite{Isbell90}.  

Factoring in Theorem~\ref{th:kc} as well (a directed terrace is a terrace),  we have~597  groups of order between 257 and 511 (inclusive) to examine.  Of these, 540 have order~384.  The program successfully worked through the others to give:  

\begin{thm}\label{th:bc}
All non-abelian groups of order $n \leq 511$, except possibly $n \in \{256, 384 \}$,  are terraced. 
\end{thm}

\noindent
Proof note.  In fact, in many cases the overhead of tracking the more complicated terrace conditions outweighed the more general nature of a terrace compared to a directed terraces and so often the directed terrace program was used. 
\qed

The result  up to order~87, except for order~64, is in~\cite{Anderson92}.

\section{Special terraces}\label{sec:special}

Recall that a {\em binary group} is a group that has a single involution.  Let G be a binary group of order~$2m$ with involution~$z$. Let~${\bf a} = (a_1, a_2, \ldots, a_{2m})$ be a directed terrace for~$G$ with associated sequencing~${\bf b} = (b_1, b_2, \ldots, b_{2m-1})$.  If $b_{m} = z$ and $b_{i} = b_{2m-i}$ for $1 \leq i \leq m-1$ then both the directed terrace and the sequencing are {\em symmetric}.  The directed terraces found by Gordon \cite{Gordon61} for abelian binary groups are symmetric.  {\em Anderson's Conjecture} is that all binary groups except the quaternion group of order~8 have symmetric sequencings.

In a series of papers \cite{Anderson87b, Anderson88, Anderson90b, AandI93, AandL88}, Anderson, Ihrig and Leonard extended Gordon's result to show that all soluble binary groups, except the quaternion group of order~8, have symmetric sequencings. 
In \cite{AandI93b} it is shown that given an insoluble binary group~$G$ with involution~$z$, the group $G/\langle z \rangle$ has a subnormal series
$$\{ e\}  \trianglelefteq G_2 \trianglelefteq G_1 \trianglelefteq G / \langle z \rangle$$
such that~$|G_2|$ and $|(G/ \langle z \rangle) / G_1|$ are odd and $G_1 / G_2$ is isomorphic to one of~$A_7$, $\psl(2,q)$ or~$\pgl(2,q)$ for odd prime powers~$q > 3$.  Call~$G_1 / G_2$ the group {\em associated} with~$G$.  It is also shown that to find a symmetric sequencing for a binary group~$G$ it is sufficient to find a terrace for the group associated with~$G$.
For each possible associated group, there are infinitely many insoluble binary groups that have it as their associated group \cite{AandI93b}.

It is already known that $\psl(2,5) \cong A_5$ is terraced \cite{nonab32}.  Theorem~\ref{th:kc} gives:

\begin{cor}
If~$G$ is an insoluble binary group whose associated group is one of
$$\psl(2,7), \psl(2,9), \psl(2,11), \pgl(2,5), \pgl(2,7) $$ 
then~$G$ has a symmetric sequencing.
\end{cor}

For undirected terraces a recent approach that allows the construction of terraces given a normal subgroup that has both even order and index in some cases is to use ``extendable terraces".

Let ${\bf a} = (a_1, a_2, \ldots, a_n)$ be a basic terrace for a group of order $n$.  Then ${\bf a}$ is {\em extendable} if $a_n = a_2^2$ and $a_{j-1}a_{j+1} = a_j = a_{j+1}a_{j-1}$ for some $j \geq 5$.  The following result does not include the full power of the machinery (which requires some fiddly notation and conditions) but gives a good sense of what is possible:

\begin{thm}\label{th:ext}{\rm \cite{Ollis15, OWext}}
Let $K$ be a group with an extendable terrace.  Suppose that~$G$  has one of the following forms:
\begin{itemize}
\item $H_1 \times H_2 \times  \cdots \times H_s \times \Z_2^t \times K$  where~$t=0$ or $t>1$ and each~$H_i \in \{\Z_6 \times \Z_2, D_8, D_{12}, A_4\} $,
\item $A \times K$ where~$A$ is an abelian $2$-group that has a normal series with each factor isomorphic to~$\Z_2^2$,  
\item $\Z_2 \times \Z_4 \times K$ and $|K|$ is congruent to $1,2,3,$ or $4 \imod{7}$.
\end{itemize}
Then~$G$ has an extendable terrace.

\end{thm}

\noindent
Proof note.  The third part of the statement is not explicitly stated elsewhere.  However, it follows immediately from Theorem~4 and the proof of Theorem~2 of~\cite{Ollis15}.
\qed


It is known that groups of order congruent to~2 modulo~4 cannot have an extendable terrace and that groups of order~6 or less or the non-abelian groups of order~8 do not have extendable terraces.  Conversely, several families of groups---including cyclic groups not ruled out in the previous sentence, dihedral groups of order a non-trivial multiple of~8, and each of the non-abelian groups of orders~12,~16 and~20---do have extendable terraces, see \cite{Ollis15}.

Rather than designing a new algorithm to attempt to look directly for extendable terraces we used a method to generate lots of terraces from a given one. 

Given a terrace, there are multiple ways to cut it into two pieces and reassemble into a new terrace (this is similar to the process of the terrace-finding algorithm, except here the goal is to maintain the perfect altitude rather than to climb).  If we forbid reversing either of the pieces but include reversing the whole terrace (which always gives a new terrace) we get an {\em orbit of terraces}. The number of essentially different terraces in an orbit must divide either~4 or~6 \cite{Bailey84}.  

By allowing moves that reverse one piece, we can form a chain of terraces that, often, encompasses many essentially different ones.  For example, of the 138,066 essentially different terraces for~$\Z_{13}$, 137,592 are mutually reachable from each other via these moves \cite{OllisPHD}.  It is perhaps better to think of the landscape in which the heuristic algorithms are operating as containing many connected ridges of the desired altitude rather than hills with distinct peaks.

The strategy for finding an extendable terrace was to  find a terrace using one of the algorithms from the previous section and then quickly generate lots of terraces for the group using these moves.  The property of being extendable is sufficiently common that we were able to reach this result:

\begin{thm}\label{th:text}
Every group of order $n$ with $12 \leq n < 96$ and $n \not\equiv 2 \imod{4}$ has an extendable terrace.
\end{thm}

In conjunction with Theorem~\ref{th:ext} (or the more general results of \cite{Ollis15}) this gives more infinite families of groups that are known to be terraced.  While this is the most effective method known to date for terracing large numbers of 2-groups, it barely scratches the surface of, for example, the~10,494,213 groups of order~512.

We now consider a stronger property than completeness that Latin squares may possess and what group theoretic structures allow us to construct such squares.

A {\em Roman-$k$ square} is a Latin square with the property that for each pair of distinct symbols $x$ and $y$ and for each integer $m$ with $1 \leq m \leq k$ we have $y$ appearing $m$ cells to the right of $x$ at most once.  A Roman-1 square is simply a Roman square.  A Roman-$(n-1)$ square is called a {\em Vatican square}.  If both a square and its transpose are Roman-$k$ then then it is called {\em $k$-complete}.

Vatican squares are known to exist only of order one less than a prime. In addition, Roman-2 squares exist for orders $2p$, where $p$ is a prime congruent to~5, 7 or~19 modulo~24, and $k=2$, and for orders $n=2m$ with $5 \leq m \leq 25$, see \cite{TuscanCRC}.  

Let ${\bf a} = (a_1, a_2, \ldots, a_n)$ be a terrace for a group~$G$ of order~$n$.  For each $m$ with $1 \leq m \leq n-1$, define ${\bf b^{(m)}} = (b_1^{(m)}, b_2^{(m)}, \ldots, b_{n-m}^{(m)})$ by $b_i^{(m)} = a_i^{-1}a_{i+m}$.  Fix $k \leq n-1$.  If for each $m \leq k$ there are no repeats within ${\bf b^{(m)}}$ then ${\bf a}$ is a {\em directed $T_k$-terrace}.  A directed $T_1$-terrace is a directed terrace and ${\bf b^{(1)}}$ is its sequencing.

If $(a_1, a_2, \ldots, a_n)$ is a  directed $T_k$-terrace for a group of order $n$ then the Latin square with $(i,j)$-entry $a_i^{-1}a_j$ is $k$-complete \cite{Anderson91a}.  Each of the known Vatican and Roman-2 squares listed above is constructed in this way from a directed $T_{k}$-terrace for a cyclic group.

As there was no known Roman-2 square of odd order, the non-abelian group of order~21---the smallest group of odd order that might admit a directed $T_2$-terrace---is of particular interest.  The heuristic algorithm was unsuccessful and so a backtracking approach was used.  Here is a directed $T_2$-terrace for $G_{21,1}$:
$$ (e, v, u, v^2, 
  u^4, u^2v, u^4v^2, u^6v^2, 
  uv, u^3, u^3v^2, u^2, 
  u^6, u^5v, u^2v^2, uv^2, 
  u^5v^2, u^4v, u^3v, u^6v, 
  u^5). $$
Hence there is a Roman-2 square of order~21.  

We were also able to find the first directed $T_2$-terraces for non-cyclic groups of even order.  The following theorem collects the results: 

\begin{thm}\label{th:t2dir}
The groups $A_4$, $Q_{12}$, $Q_{16}$, $G_{16,6}$, $G_{16,13}$  and all non-abelian groups of orders~$n$ with $18 \leq n \leq 22$ have a directed $T_2$-terrace.  Other non-abelian groups of order at most~$22$ do not have a directed $T_2$-terrace.  Also,  $\Z_6 \times \Z_3$ and $S_4$ have directed $T_2$-terraces.  

The only groups of order less than~$20$ that have a directed~$T_3$-terrace are the cyclic groups of order~$p-1$ for prime~$p$.
\end{thm}

\noindent
Proof note.  The terraces that prove this are available at the website associated with this paper, the link for which is given in the introduction.  As an example, here is the directed $T_2$-terrace for~$A_4$:
\begin{flushleft}
$( (), (2,3,4), (1,2)(3,4), (1,3)(2,4), (1,3,4), (1,4,3),$
\end{flushleft}
\begin{flushright}
$ (1,2,4), 
  (1,4,2), (1,3,2), (1,2,3), (2,4,3), (1,4)(2,3) ).$ \qed
\end{flushright}

The final type of terrace that has attracted interest that we consider here is the half-and-half terrace and the stronger narcissistic terrace introduced in \cite{AP02}. 

 Let ${\bf a} = (a_1, a_2, \ldots, a_n)$ be a terrace for a group~$G$ of odd order~$n$ and ${\bf b}= (b_1, b_2, \ldots, b_{n-1})$ be its 2-sequencing.  If there is one occurrence of from each set  $\{ g, g^{-1} :  g^2 \neq e \}$ in $(b_1, b_2, \ldots, b_{(n-1)/2})$, and hence also one occurrence from that set in $(b_{(n+1)/2}, \ldots, b_n)$, then ${\bf a}$ is a {\em half-and-half terrace}.   If ${\bf b}$ is equal to its reverse then ${\bf b}$ is {\em reflective} and ${\bf a}$ is {\em narcissistic}.  Narcissistic terraces are necessarily half-and-half.

It is known that all abelian groups of odd order have a narcissistic terrace and that $G_{21,1}$ has a directed half-and-half terrace and a narcissistic terrace~\cite{OS05}.   For non-abelian groups, the following result enables the construction of infinite families of examples of half-and-half and narcissistic terraces:

\begin{thm}\label{th:hh}{\rm \cite{OS05}}
Let~$G$ and~$H$ be groups of odd order.
If~$G$ and~$H$ have half-and-half terraces then~$G \times H$ has a half-and-half terrace. If~$G$ and~$H$ have narcissistic terraces then~$G \times H$ has a narcissistic terrace. \end{thm}

Similarly to the case with directed $T_2$-terraces, a backtracking approach was more successful than the heuristic algorithm.

\begin{thm}\label{th:hh}
Each of the two non-abelian groups of order~$27$ has both a narcissistic and a directed half-and-half terrace.  The non-abelian group~$G_{39,1}$ of order~$39$ has a directed half-and-half terrace.
\end{thm}

\noindent
Proof note.  The terraces that prove the result are available at the website associated with this paper, the link for which is given in the introduction.    We include here a narcissistic terrace and a directed half-and-half terrace for~$G_{27,4}$:

\begin{flushleft}
$e, v, u, v^2, 
  u^4v^2, u^5v, u^4, 
  u^5v^2, u^6v^2, u^3v^2, uv^2, 
  u^8, u^2v, u^5, $
  \end{flushleft}
  \begin{flushright}
$  u^8v^2, u^2, u^6v, 
  u^7v, u^4v, u^8v, 
  u^6, u^2v^2, u^3v, 
  uv, u^3, u^7v^2, 
  u^7)$ 
  \end{flushright}

\begin{flushleft}
$e, v, u, v^2, 
  u^4v^2, u^7v, u^8v^2, 
  u^2, u^3v^2, u^5v, 
  uv, u^5, u^8,$
  \end{flushleft}
  \begin{flushright} 
$  u^3, u^8v, u^4, 
  uv^2, u^6v^2, u^7, 
  u^4v, u^2v^2, u^2v, 
  u^6v, u^3v, u^5v^2, u^7v^2, u^6$
\qed
\end{flushright}

\section{Very small groups}\label{sec:vsg}

In this section we report on enumerations of terraces and directed terraces in groups of order up to~15.  Let~$t(G)$ denoted the number of essentially different terraces of~$G$ and $d(G)$ denote the number of essentially different directed terraces of~$G$.

For groups of order at most 9, $\Z_{10}$, $\Z_{11}$ and~$\Z_{13}$ the values of $t(G)$ and $d(G)$ have already been computed, sometimes by hand \cite{Bailey84, OllisPHD, Street88}.   There are two known errors in these computations: \cite{Street88} gives $t(\Z_5) = 2$ rather than~3 and \cite{Bailey84} gives $t(\Z_3^2) = 32$ rather than 35.  Other than those, the new computation matched these results.  We give the results of the  computation for $5 \leq |G| \leq 15$ in Table~\ref{tab:enum}.    Each of~$\Z_2$, $\Z_3$ and~$\Z_4$ has exactly one terrace, which is directed for~$\Z_2$ and~$\Z_4$.  Elementary abelian 2-groups have no terraces.

\begin{table}
\caption{Enumerations for $5 \leq |G| \leq 15$}\label{tab:enum}
$$
\begin{array}{rrr}
\hline
G & t(G) & d(G) \\
\hline
\Z_5 & 3 & 0 \\
\Z_6 &  11 & 2 \\
D_6 & 2 & 0 \\  
\Z_8 & 58 & 6 \\
\Z_4 \times Z_2 & 10 & 0 \\
D_8  & 6 &  0 \\
Q_8 & 6 &  0 \\
\Z_9 & 234 & 0 \\
\Z_3 \times \Z_3 & 35 & 0 \\
\Z_{10} & 1,517 & 72 \\
D_{10} & 76 & 16 \\
\Z_{11} &  4,116 & 0 \\
\Z_{12} &  40,722 & 964 \\
\Z_6 \times \Z_2 &  5,528 & 0 \\
D_{12} & 1,380 & 256 \\
Q_{12} & 13,470 & 372 \\
A_4 & 3,516 & 96 \\
\Z_{13} & 138,066 &  0 \\
\Z_{14} & 1,458,038 & 14,888 \\
D_{14} & 25,608 & 2,700 \\
\Z_{15} & 10,910,262 & 0 \\
\hline
\end{array}
$$
\end{table}

\end{document}